\documentclass[11pt]{article} % use larger type; default would be 10pt
\usepackage{mathrsfs}
\usepackage[flushleft]{threeparttable}
\usepackage{makecell,booktabs}
\usepackage{multirow}
\usepackage[utf8]{inputenc}
\usepackage{algorithm}
\usepackage[noend]{algpseudocode}

\usepackage{geometry} % to change the page dimensions
\usepackage[utf8]{inputenc}
\usepackage{mdframed} % cool box 
\usepackage{graphicx} % support the \includegraphics command and options
\usepackage{booktabs} % for much better looking tables
\usepackage{array} % for better arrays (eg matrices) in maths
\usepackage{paralist} % very flexible & customisable lists (eg. enumerate/itemize, etc.)
\usepackage{verbatim} % adds environment for commenting out blocks of text & for better verbatim
\usepackage{subfig} % make it possible to include more than one captioned /table in a single float
\usepackage{amssymb, amsmath, amsthm,tikz,listings,esint,cancel,mathabx, mathtools}

\usepackage{capt-of}
\usepackage{setspace}
\usepackage{hyperref, comment}

\usepackage[
backend=biber,
style=numeric,
sorting=nyt
]{biblatex}
\renewcommand{\emptyset}{\font\cmsy = cmsy10 at 10pt
 \hbox{\cmsy \char 59}
}

\renewcommand{\emptyset}{\font\cmsy = cmsy10 at 10pt
 \hbox{\cmsy \char 59}
}

\usepackage{sectsty}
\allsectionsfont{\sffamily\mdseries\upshape} % (See the fntguide.pdf for font help)
\usepackage[nottoc,notlof,notlot]{tocbibind} % Put the bibliography in the ToC

\begin{document}
\begin{titlepage}
\pagenumbering{roman}
   \begin{center}
       %\vspace*{1cm}
        \Huge
       \textbf{On the Truth Assignment Theorem of the Language of Sentential Logic}
\Large

    \vspace{3cm}
       
       \textbf{A CAPSTONE PROJECT SUBMITTED TO THE FACULTY OF THE SCHOOL OF MATHEMATICS OF THE UNIVERSITY OF MINNESOTA BY}

       \vspace{3cm}
       
       \textbf{Tianyi Sun}

        \vspace{3cm}

    \textbf{ IN PARTIAL FULFILLMENT OF THE REQUIREMENTS FOR THE DEGREE OF \\ BACHELOR  OF ARTS IN MATHEMATICS}
    
    \vspace{2cm}
       \textbf{Karel Prikry}
      
       \vspace{2cm}
       
      % \textbf{University of Minnesota Twin - Cities\\}
       \textbf{July 13, 2021}
         %\vspace{1.5cm}

   \end{center}
\end{titlepage}

\begin{titlepage}
\pagenumbering{roman}
   \begin{center}
%\begin{copyrightpage }
    \vfill

   \vspace*{9cm}
   
\textbf{\textcopyright Tianyi Sun 2021 \\ ALL RIGHTS RESERVED}

   \end{center}
\end{titlepage}

\tableofcontents
\pagenumbering{arabic}
 \clearpage

\section{Introduction}
This project contains two chapters. Chapter 2 has two sections. First, we define the well-formed formulas of the Language of Sentential Logic using Construction Sequences. Second, we prove the Truth Assignments using the Language of Sentential Logic. Chapter 3 has two sections. First, we define the Recursion Theorem. Second, we prove the Truth Assignments using a general formulation of the Recursion Theorem.

\section{Construction Sequences}
We want to go through the proof of the Truth Assignment Theorem which is the Theorem 12A in \cite{enderton2001mathematical}, and we will also discuss some related results. 

We will first work out a rigorous definition of well-formed formulas (or simply formulas or wffs) of Language of Sentential Logic (LSL), using the Construction Sequences, as mentioned on pages 17 and 35 in \cite{enderton2001mathematical}. 

Firstly, we define the operations $\varepsilon_\neg$, $\varepsilon_\wedge$, $\varepsilon_\vee$, $\varepsilon_\rightarrow$, and
$\varepsilon_\leftrightarrow$, which are operations or arbitrary LSL expressions, by proceeding exactly as at the top of page 17 in \cite{enderton2001mathematical}: If $\alpha$ and $\beta$ are arbitrary expressions (a finite sequence of symbols) of LSL (i.e. not just wffs), then 
\begin{align*}
    \varepsilon_\neg(\alpha) &= (\neg \alpha), \\
    \varepsilon_\wedge(\alpha, \beta) &= (\alpha \wedge \beta),\\
    \varepsilon_\vee(\alpha, \beta) &= (\alpha \vee \beta),\\
    \varepsilon_\rightarrow(\alpha, \beta) &= (\alpha \rightarrow \beta),\\
    \varepsilon_\leftrightarrow(\alpha, \beta) &= (\alpha \leftrightarrow \beta).
\end{align*}

\smallskip
\noindent{\bf Construction Sequences.}
The definition of Construction Sequences is on pages 17 and 18 in \cite{enderton2001mathematical}. 
A Construction Sequence is a finite sequence $\langle\varepsilon_1, \varepsilon_2, ..., \varepsilon_n \rangle$ of expressions, such that for every $i\preceq n$, $\varepsilon_i$ is either a sentence symbol, or there are some $j \prec i$ such that $\varepsilon_i = \varepsilon_\neg(\varepsilon_j)$, or there are $j\prec i$ and $k \prec i$ such that $\varepsilon_i = \varepsilon_\square(\varepsilon_j, \varepsilon_k)$, where $\square$ is one of $\neg, \wedge, \vee, \rightarrow,$ and $\leftrightarrow$. 

We now define wffs as follow. 

\smallskip
\noindent{\bf Well-formed formulas (wffs).} An expression $\alpha$ is a wff, if there exists a construction sequence $\langle\varepsilon_1, \varepsilon_2, ..., \varepsilon_n \rangle$ such that $\alpha = \varepsilon_n$.

As in \cite{enderton2001mathematical}, if $S$ is any set of sentence symbols, then $\overline{S}$ is the set of all wffs $\alpha$ such that every sentence symbol occurring in $\alpha$ belongs to $S$. Also, $\alpha \in \overline{S}$, if and only if $\alpha$ has a construction sequence $\langle\varepsilon_1, \varepsilon_2,..., \varepsilon_n\rangle$, with $\varepsilon_n = \alpha$, and every sentence symbol, belonging to this construction sequence, belongs to $S$. (See also in \cite{enderton2001mathematical}, page 20, the paragraph above the Condition 0.) We will call such construction sequences S-based. 

We will now consider an arbitrary truth assignment $$v:S \rightarrow \{T,F\}$$ and we will define, for every S-based construction sequence $\langle\varepsilon_1, \varepsilon_2,..., \varepsilon_n\rangle$, an associated sequence of truth values $\langle\tau_1, \tau_2,..., \tau_n\rangle$. We will use induction on positive integers $i\preceq n$. Thus, for $i=1$, $\varepsilon_1$ has to be a sentence symbol belonging to $S$, and we will set $\tau_1 = v(\varepsilon_1)$.

Now suppose that $i\preceq n$, and $\tau_j$ is defined for every $j \prec i$. We now consider all cases as in the definition of the construction sequence: If $\varepsilon_i$ is a sentence symbol belonging to $S$, we set $\tau_i = v(\varepsilon_i)$, exactly as when $i=1$, following the condition 0., at the bottom of page 20, in \cite{enderton2001mathematical}. If $\varepsilon_i = (\neg \varepsilon_j)$ (i.e. $\varepsilon_\neg(\varepsilon_j)$) for some $j \prec i$, we set 

$$\tau_i =
\setlength\arraycolsep{1pt}
\left\{
\begin{array}{rcrcrc@{\qquad}l}
F &\mbox{ if } \tau_j = T, \\ 
T &\mbox{ if } \tau_j = F.\\ 
\end{array}
\right.
$$
If there exist $j \prec i$ and $k \prec i$ such that $\varepsilon_i = \varepsilon_\square(\varepsilon_j, \varepsilon_k)$ , where $\square \in \{ \neg, \wedge,\vee, \rightarrow, \leftrightarrow \}$, then $\tau_i$ is obtained from $\tau_j, \tau_k$ according to the corresponding condition from among the conditions 2.- 5. at the top of page 21, in \cite{enderton2001mathematical}. For example, if $\square = \rightarrow$, then following the condition 4., we set

$$\tau_i =
\setlength\arraycolsep{1pt}
\left\{
\begin{array}{rcrcrc@{\qquad}l}
F &\mbox{ if } \tau_j = T \mbox{ and } \tau_k = F,\\ 
T &\mbox{ otherwise.}\\ 
\end{array}
\right.
$$

Thus, the method of induction allows us to conclude that $\tau_i$ is defined for all $i\preceq n$. See also the following example. 

{\bf Example.} Let the construction sequence be $\langle A_1, A_2, A_3, (\neg A_2), (A_3 \vee (\neg A_2)), ((A_3 \vee (\neg A_2))\leftrightarrow A_1)\rangle$. Let $v: \{A_1,A_2,A_3\} \rightarrow \{F,T,F\}$ be $v(A_1) = F, v(A_2) = T, v(A_3) = F$. Then the associated sequence of truth values is $\tau_1 = F, \tau_2 = T, \tau_3 = F$.
Then $\tau_4 = F, \tau_5 = F, \tau_6 = T$. 
Hence, $\langle\tau_1, \tau_2, \tau_3, \tau_4, \tau_5, \tau_6\rangle = \langle F, T, F, F, F, T\rangle$.

\smallskip
\noindent{\bf The Main Lemma.} Let $\alpha \in \overline{S}$ and $\langle\delta_1,\delta _2,...,\delta_m\rangle$ and $\langle\varepsilon_1, \varepsilon_2,..., \varepsilon_n\rangle$ be S-based construction sequences for $\alpha$, i.e. $\alpha = \delta_m = \varepsilon_n $. Furthermore let $\langle\tau_1, \tau_2,...,\tau_m\rangle$ and $\langle\tau_1', \tau_2',...,\tau_n'\rangle$ be the associated sequences of truth values. Then $\tau_m = \tau_n'$. 

The proof uses the Induction Principle (page 18, in \cite{enderton2001mathematical}), and also the Unique Readability Theorem as it is stated in the MATH5165 class notes, and also on p.40 in \cite{enderton2001mathematical}.

\smallskip
\noindent{\bf Unique Readability Theorem.} Let $\alpha$ be a wff of the language sentential logic. Then exactly one of the following six possibilities happens: 

%\begin{multicols}
 \begin{enumerate}
    \item  there is a unique natural number $n$, such that $\alpha$ is the sentence symbol $A_n$; 
    \item  there is a unique wff $\beta$, such that $\alpha$ is $(\neg \beta)$;
    \item  there are unique wffs $\beta, \gamma$, such that $\alpha$ is $(\beta \wedge \gamma)$;
    \item  there are unique wffs $\beta, \gamma$, such that $\alpha$ is $(\beta \vee \gamma)$; 
    \item  there are unique wffs $\beta, \gamma$, such that $\alpha$ is $(\beta \rightarrow \gamma)$;
    \item  there are unique wffs $\beta, \gamma$, such that $\alpha$ is $(\beta \leftrightarrow \gamma)$.
 \end{enumerate}   
%\end{multicols}

{\bf Comment.} The Unique Readability Theorem is equivalent to saying that every wff has exactly one parsing by matching appropriate pairs of left and right parentheses. Or, using the approach in the Book, every wff has exactly one ancestral tree. However, as far as a precise mathematical statement is concerned, we will stay with the formulation of the Unique Readability Theorem given above. The reason for this is that otherwise we would need to provide a more precise definition of parsing.

\smallskip
\noindent{\bf Proof of the Main Lemma.}

We want to state the Induction Principle in slightly stronger form, as needed for the proof of the Lemma above, and of the Theorem 12A. Namely instead of just dealing with the set of all wffs, we need to state the Induction Principle for the set $\overline{S}$, where $S$ is any set of sentence symbols. 

Thus, 
\smallskip
\noindent{\bf the Induction Principle for $S$:} 

Let $S$ be any set of sentence symbols and $\mathscr{S} \subseteq \overline{S}$ be such that 
\begin{enumerate}
\item $S \subseteq \mathscr{S}$, and

\item $\mathscr{S}$ is closed under the five formula building operations. 

Then $\mathscr{S} = \overline{S}$.
\end{enumerate}
We set up the proof using the Induction Principle by defining $\mathscr{S}$ to be the set of those wff $\alpha \in \overline{S}$ for which the statement of the Lemma holds. Hence, we need to prove that
\begin{enumerate}
        \item  $S \subseteq \mathscr{S}$, and
        \item  $\mathscr{S}$ is closed under the five formula building operations. 
\end{enumerate}

The Induction Principle then allows us to conclude that $\mathscr{S} = \overline{S}$. Thus, the statement of the Lemma holds for every $\alpha \in \overline{S}$, as designed. 
\begin{proof}

Proof of 1.. Let $\langle\delta_1,\delta _2,...,\delta_m\rangle$ and $\langle\varepsilon_1, \varepsilon_2,..., \varepsilon_n\rangle$ be construction sequences for a sentence symbol $A_j$ belonging to $S$. Thus $\delta_m = \varepsilon_n = A_j$. But then by the definition of the associated construction sequences, we obtain $\tau_m = v(A_j)$ and $\tau_n^{'} = v(A_j)$, hence $\tau_m = \tau_n^{'}$. 

Proof of 2.. We first consider domain under $\neg$. Hence, let $\alpha \in \mathscr{S}$. We want to prove that $(\neg \alpha) \in \mathscr{S}$. Thus, let $\langle \delta_1,\delta _2,...,\delta_m\rangle$ and $\langle\varepsilon_1, \varepsilon_2,..., \varepsilon_n\rangle$ be construction sequences for $(\neg \alpha)$. Thus $\delta_m = \varepsilon_n = (\neg \alpha)$, and there exist $j\langle m$ and $k\langle n$ such that $\delta_j = \alpha$ and $\varepsilon_k = \alpha$. Hence, $\langle\delta_1,\delta_2,...,\delta_j\rangle$ and $\langle\varepsilon_1, \varepsilon_2,..., \varepsilon_k\rangle$ both are construction sequences for $\alpha$, and $\langle\tau_1,\tau_2,...,\tau_j\rangle$ and $\langle\tau_1^{'}, \tau_2^{'},..., \tau_k^{'}\rangle$ are the corresponding associated sequences of truth values. Thus, since $\alpha \in \mathscr{S}$, we obtain $\tau_j = \tau_k^{'}$. But $\delta_m = (\neg \delta_j)$ and $\varepsilon_n = (\neg \varepsilon_k)$. Hence, if $\tau_j = \tau_k^{'} =T$, we obtain $\tau_m = \tau_n^{'} =F$, and similarly for $\tau_j = \tau_k^{'} =F$. Thus, $\tau_m = \tau_n^{'}$. Therefore, $(\neg \delta) \in \mathscr{S}$.

It remains to be shown that $\mathscr{S}$ is closed under the formula building operations corresponding to the binary sentential connectives. We will do in detail the case of $\wedge$. The other three cases are similar, so we will omit them. 

Thus, let $\alpha, \beta$ both belong to $\mathscr{S}$. We want to show that $(\alpha \wedge \beta)$ belongs to $\mathscr{S}$. Thus, let $\langle \delta_1,\delta_2,...,\delta_m \rangle$ and $\langle \varepsilon_1, \varepsilon_2,..., \varepsilon_n \rangle$ be construction sequences for $(\alpha \wedge \beta)$. Hence, $\delta_m = \varepsilon_n = (\alpha \wedge \beta)$ and there exist $i\prec m$ and $j\prec m$, and $k\prec n$ and $l\prec n$, such that $\delta_i = \alpha$ and $\delta_j = \beta$, as well as $\varepsilon_k = \alpha$ and $\varepsilon_l = \beta$. But then, similarly as with the case of $\neg$, we conclude that $\tau_i = \tau_k'$ and $\tau_j = \tau_l'$, since $\alpha$ and $\beta$ both belong to $\mathscr{S}$. But then $\tau_m = \tau_n'$, since $\tau_m$ is obtained from $\tau_i$ and $\tau_j$ so that the condition 3. at the top of page 21 in \cite{enderton2001mathematical} would hold, and $\tau_n'$ is obtained in a similar manner from $\tau_k'$ and $\tau_l'$. Hence, $(\alpha \wedge \beta) \in \mathscr{S}$.

Thus, as already stated at the beginning of the proof, we conclude that $\mathscr{S} = \overline{S}$, and the Lemma is proved. 

\end{proof}
Proceeding to the next stage of the proof of Theorem 12A, we define a set $\mathscr{P}$ of ordered pairs $\langle\alpha,V\rangle$, where $\alpha \in \overline{S}$ and $V \in \{F,T\}$ as follows: $\langle\alpha,V\rangle$ belongs to $\mathscr{P}$ if and only if there exists a S-based construction sequence $\langle\delta_1,\delta_2,...,\delta_m\rangle$ for $\alpha$ such that the associated sequence of truth values $\langle\tau_1,...,\tau_m\rangle$ has $\tau_m = V$. 

{\bf Comment.} Thus for example, if $\alpha$ is $((A_3 \vee (\neg A_2)) \leftrightarrow A_1)$, then the Example preceding the statement of the Main Lemma, shows that $\langle\alpha,T\rangle \in \mathscr{P}$, where $\alpha = ((A_3 \vee(\neg A_2))\leftrightarrow A_1)$, since the last entry, $\tau_6$, of the associated sequence of truth values, is $T$. 

{\bf The Uniqueness Lemma.} For every $\alpha$ belonging to $\overline{S}$ there exists exactly one $V\in\{F,T\}$ such that $\langle\alpha,V\rangle \in \mathscr{P}$. 

\begin{proof}

Let $\alpha$, belonging to $\overline{S}$, be arbitrary. Thus by definition of wffs, $\alpha$ has a construction sequence, let us say $\langle\delta_1, \delta_2, ..., \delta_m\rangle$, and we can assume that this construction sequence is S-based. Let $\langle\tau_1,\tau_2,...,\tau_m\rangle$ be the associated sequence of truth values. Thus $\langle\alpha,V\rangle$, where $V=\tau_m$, belongs to $\mathscr{P}$. Moreover, suppose that $\langle\alpha,V\rangle$ and $\langle\alpha,V'\rangle$ both belong to $\mathscr{P}$. Thus there exist S-based construction sequences, $\langle\delta_1,\delta_2,...,\delta_m\rangle$ and $\langle\varepsilon_1,\varepsilon_2,...,\varepsilon_n\rangle$, for $\alpha$, and the associated sequences, $\langle\tau_1,\tau_2,...,\tau_m\rangle$ and $\langle\tau_1',\tau_2',...,\tau_n'\rangle$, of truth values, such that $\tau_m = V$ and $\tau_n' = V'$. But then, we obtain, by the Main Lemma, that $V = V'$. Hence the Lemma is proved. 

\end{proof}
%2.6-------
The Uniqueness Lemma shows that the set $\mathscr{P}$ satisfies the conditions for defining a function with domain $\overline{S}$ and co-domain $\{T,F\}$, or, formally, being a subset of $\overline{S} \times \{T,F\}$, it already is such a function. However, for clarity, we will write down the definition of this function, calling it $h$. Thus, for an arbitrary $\alpha$ in $\overline{S}$, and $V \in \{T,F\}$, we will define $$h(\alpha) = V,$$ if and only if $\langle\alpha, V\rangle \in \mathscr{P}$. 

Thus clearly, the Uniqueness Lemma is precisely what is needed for $h$ to be a function with domain $\overline{S}$ and co-domain $\{T,F\}$. 

We will now prove that the function $h$ satisfies the conditions 0.-5. on pages 20 and 21 in \cite{enderton2001mathematical}, that are required of $\overline{v}$. Thus a function satisfying the conditions 0.-5., required of $\overline{v}$, exists, and we can set $\overline{v} = h$. Thus in order to complete the existence part of the proof of the Truth Assignment Theorem, i.e. the Theorem 12A in \cite{enderton2001mathematical}, it remains to show that the conditions 0.-5. hold for $h$. We now proceed to do so. 

\begin{proof}

\underline{Condition 0.} Suppose that $\alpha$ is a sentence symbol belonging to $S$, i.e. $\alpha = A_n \in S$. Hence setting $\varepsilon_1 = A_n$, we obtain the construction sequence $\langle\varepsilon_1\rangle = \langle A_n\rangle$ for $\alpha = A_n$, and its associated sequence of truth values is $\langle \tau_1 \rangle = \langle v (A_n)\rangle$. Thus $\langle\varepsilon_1,\tau_1\rangle = \langle A_n,v(A_n)\rangle$ belongs to $\mathscr{P}$. Therefore $h(A_n) = v(A_n)$, which shows that the Condition 0. holds for $h$. 

\underline{Condition 1.} Let $\alpha \in \overline{S}$. We need to show that 
$$h(( \neg \alpha)) =
\setlength\arraycolsep{1pt}
\left\{
\begin{array}{rcrcrc@{\qquad}l}
T &\mbox{ if } h(\alpha) = F\\ 
F &\mbox{ if } h(\alpha) = T.\\ 
\end{array}
\right.
$$
Thus suppose $h(\alpha) = F$. Hence $\langle\alpha, F\rangle \in \mathscr{P}$, and therefore there exists a construction sequence $\langle\delta_1,\delta_2,...,\delta_m\rangle$ for $\alpha$ and its associated sequence of truth values $\langle\tau_1,\tau_2,...,\tau_m\rangle$ with the properties $\delta_m = \alpha$ and $\tau_m = F$. But then $\langle\delta_1,\delta_2,...,\delta_m, ( \neg \alpha)\rangle$ is a construction sequence for $(\neg \alpha)$, and its associated sequence of truth values is $\langle \tau_1,\tau_2,...,\tau_m, T\rangle$, since $\tau_m = F$. Hence, by the definition of $\mathscr{P}$, we obtain $\langle(\neg \alpha), T\rangle \in \mathscr{P}$. Hence, $$h((\neg \alpha)) = T,$$ as desired. 

We similarly prove that if $h(\alpha) = T$, then $h((\neg \alpha)) = F$. 

Hence the Condition 1. holds for $h$. 

\underline{Condition 2.} Suppose that $\alpha$ and $\beta$ both belong to $\overline{S}$. We need to show that 
$$h((\alpha \wedge \beta)) =
\setlength\arraycolsep{1pt}
\left\{
\begin{array}{rcrcrc@{\qquad}l}
T &\mbox{ if } h(\alpha) = h(\beta) = T\\ 
F & \mbox{ otherwise }.\\ 
\end{array}
\right.
$$

We first suppose that $h(\alpha) = h(\beta) = T$. We need to prove that $$h((\alpha \wedge \beta)) = T.$$
Hence $\langle\alpha,T\rangle$ and $\langle\beta,T\rangle$
both belong to $\mathscr{P}$, and therefore there exists a construction sequence $\langle\delta_1,\delta_2,...,\delta_m\rangle$ for $\alpha$ and its associated sequence of truth values $\langle\tau_1,\tau_2,...,\tau_m\rangle$ with the properties $\delta_m = \alpha$ and $\tau_m = T$. And there exists another construction sequence $\langle\varepsilon_1,\varepsilon_2, ...,  \varepsilon_n\rangle$ for $\beta$ and its associated sequence of truth values $\langle\tau_1^{'},\tau_2^{'}, ..., \tau_n^{'}\rangle$ with the properties $\varepsilon_n = \beta$ and $\tau_n^{'} = T$. Hence $\langle \delta_1,\delta_2,...,\delta_m,\varepsilon_1,\varepsilon_2, ..., \varepsilon_n, (\alpha \wedge \beta) \rangle$ is a construction sequence for $(\alpha \wedge \beta)$. 
Because the condition 2. holds for the associated truth value sequences, the associated sequence of truth values for $\langle\delta_1,\delta_2,...,\delta_m,\varepsilon_1,\varepsilon_2, ..., \varepsilon_n, (\alpha \wedge \beta) \rangle$ is
$\langle\tau_1,\tau_2,...,T,\tau_1^{'},\tau_2^{'}, ..., T, T\rangle $, since $\tau_m = T$ and $\tau_n^{'} = T$. 
Hence, by the definition of $\mathscr{P}$, we obtain $\langle(\alpha \wedge \beta), T\rangle \in \mathscr{P}$. Hence $h((\alpha \wedge \beta)) = T$, as desired. 

We second suppose that $h(\alpha) = T$ and $h(\beta) = F$. We need to prove that $$h((\alpha \wedge \beta)) = F.$$
Hence $\langle\alpha,T\rangle$ and $\langle\beta,F\rangle$
both belong to $\mathscr{P}$, and therefore there exists a construction sequence $\langle\delta_1,\delta_2,...,\delta_m\rangle$ for $\alpha$ and its associated sequence of truth values $\langle\tau_1,\tau_2,...,\tau_m\rangle$ with the properties $\delta_m = \alpha$ and $\tau_m = T$. And there exists another construction sequence $\langle\varepsilon_1,\varepsilon_2, ..., \varepsilon_n\rangle$ for $\beta$ and its associated sequence of truth values $\langle\tau_1^{'},\tau_2^{'}, ..., \tau_n^{'}\rangle$ with the properties $\varepsilon_n = \beta$ and $\tau_n^{'} = F$. Hence $\langle\delta_1,\delta_2,...,\delta_m,\varepsilon_1,\varepsilon_2, ..., \varepsilon_n, (\alpha \wedge \beta) \rangle$ is a construction sequence for $(\alpha \wedge \beta)$. 
Because the condition 2. holds for the associated truth value sequences, 
the associated sequence of truth values for $\langle\delta_1,\delta_2,...,\delta_m,\varepsilon_1,\varepsilon_2, ..., \varepsilon_n, (\alpha \wedge \beta) \rangle$ is
$\langle\tau_1,\tau_2,...,T,\tau_1^{'},\tau_2^{'}, ..., F, F\rangle $,
since $\tau_m = T$ and $\tau_n^{'} = F$.
Hence, by the definition of $\mathscr{P}$, we obtain $\langle(\alpha \wedge \beta), F\rangle \in \mathscr{P}$. Hence $h((\alpha \wedge \beta)) = F$, as desired. 

We similarly prove that if $h(\alpha) = F$ and $h(\beta) = T$, then $$h((\alpha \wedge \beta)) = F.$$

We third suppose that $h(\alpha) = F$ and $h(\beta) = F$. We need to prove that $$h((\alpha \wedge \beta)) = F.$$
Hence $\langle\alpha,F\rangle$ and $\langle\beta,F\rangle$
both belong to $\mathscr{P}$, and therefore there exists a construction sequence $\langle\delta_1,\delta_2,...,\delta_m\rangle$ for $\alpha$ and its associated sequence of truth values $\langle\tau_1,\tau_2,...,\tau_m\rangle$ with the properties $\delta_m = \alpha$ and $\tau_m = F$. And there exists another construction sequence $\langle\varepsilon_1,\varepsilon_2, ..., \varepsilon_n\rangle$ for $\beta$ and its associated sequence of truth values $\langle\tau_1^{'},\tau_2^{'}, ..., \tau_n^{'}\rangle$ with the properties $\varepsilon_n = \beta$ and $\tau_n^{'} = F$. Hence $\langle\delta_1,\delta_2,...,\delta_m,\varepsilon_1,\varepsilon_2, ..., \varepsilon_n, (\alpha \wedge \beta) \rangle$ is a construction sequence for $(\alpha \wedge \beta)$. Because the condition 2. holds for the associated truth value sequences, the associated sequence of truth values for $\langle\delta_1,\delta_2,...,\delta_m,\varepsilon_1,\varepsilon_2, ..., \varepsilon_n, (\alpha \wedge \beta) \rangle$ is $\langle\tau_1,\tau_2,...,F,\tau_1^{'},\tau_2^{'}, ..., F, F\rangle $, since $\tau_m = F$ and $\tau_n^{'} = F$. Hence, by the definition of $\mathscr{P}$, we obtain $\langle(\alpha \wedge \beta), F\rangle \in \mathscr{P}$. Hence $h((\alpha \wedge \beta)) = F$, as desired.

Therefore, Condition 2. holds for $h$.

\underline{Condition 3.} Suppose that $\alpha$ and $\beta$ both belong to $\overline{S}$. We need to show that 
$$h((\alpha \vee \beta)) =
\setlength\arraycolsep{1pt}
\left\{
\begin{array}{rcrcrc@{\qquad}l}
T &\mbox{ if } h(\alpha) = T \mbox{ or } h(\beta) = T \mbox{ (or both) }\\ 
F & \mbox{ otherwise }.\\ 
\end{array}
\right.
$$

We first suppose that $h(\alpha) = T$. We need to prove that if either $h(\beta) = T$ or $h(\beta) = F$, we always have $$h((\alpha \vee \beta)) = T.$$ Thus, we prove this by cases. 

Case1. Assume $h(\beta) = T$. Hence $\langle\alpha,T\rangle$ and $\langle\beta,T\rangle$
both belong to $\mathscr{P}$, and therefore there exists a construction sequence $\langle\delta_1,\delta_2,...,\delta_m\rangle$ for $\alpha$ and its associated sequence of truth values $\langle\tau_1,\tau_2,...,\tau_m\rangle$ with the properties $\delta_m = \alpha$ and $\tau_m = T$. And there exists another construction sequence $\langle\varepsilon_1,\varepsilon_2, ...,  \varepsilon_n\rangle$ for $\beta$ and its associated sequence of truth values $\langle\tau_1^{'},\tau_2^{'}, ..., \tau_n^{'}\rangle$ with the properties $\varepsilon_n = \beta$ and $\tau_n^{'} = T$. Hence $\langle \delta_1,\delta_2,...,\delta_m,\varepsilon_1,\varepsilon_2, ..., \varepsilon_n, (\alpha \vee \beta) \rangle$ is a construction sequence for $(\alpha \vee \beta)$. 
Because the condition 3. holds for the associated truth value sequences, the associated sequence of truth values for $\langle\delta_1,\delta_2,...,\delta_m,\varepsilon_1,\varepsilon_2, ..., \varepsilon_n, (\alpha \vee \beta) \rangle$ is
$\langle\tau_1,\tau_2,...,T,\tau_1^{'},\tau_2^{'}, ..., T, T\rangle $, since $\tau_m = T$ and $\tau_n^{'} = T$. 
Hence, by the definition of $\mathscr{P}$, we obtain $\langle(\alpha \vee \beta), T\rangle \in \mathscr{P}$. 

Case2. Assume $h(\beta) = F$. Hence $\langle\alpha,T\rangle$ and $\langle\beta,F\rangle$
both belong to $\mathscr{P}$, and therefore there exists a construction sequence $\langle\delta_1,\delta_2,...,\delta_m\rangle$ for $\alpha$ and its associated sequence of truth values $\langle\tau_1,\tau_2,...,\tau_m\rangle$ with the properties $\delta_m = \alpha$ and $\tau_m = T$. And there exists another construction sequence $\langle\varepsilon_1,\varepsilon_2, ...,  \varepsilon_n\rangle$ for $\beta$ and its associated sequence of truth values $\langle\tau_1^{'},\tau_2^{'}, ..., \tau_n^{'}\rangle$ with the properties $\varepsilon_n = \beta$ and $\tau_n^{'} = F$. Hence $\langle \delta_1,\delta_2,...,\delta_m,\varepsilon_1,\varepsilon_2, ..., \varepsilon_n, (\alpha \vee \beta) \rangle$ is a construction sequence for $(\alpha \vee \beta)$. 
Because the condition 3. holds for the associated truth value sequences, the associated sequence of truth values for $\langle\delta_1,\delta_2,...,\delta_m,\varepsilon_1,\varepsilon_2, ..., \varepsilon_n, (\alpha \vee \beta) \rangle$ is
$\langle\tau_1,\tau_2,...,T,\tau_1^{'},\tau_2^{'}, ..., F, T\rangle $, since $\tau_m = T$ and $\tau_n^{'} = F$. 
Hence, by the definition of $\mathscr{P}$, we obtain $\langle(\alpha \vee \beta), T\rangle \in \mathscr{P}$. 

Hence $h((\alpha \vee \beta)) = T$, as desired. 

We similarly prove that if $h(\beta) = T$ and either $h(\alpha) = T$ or $h(\alpha) = F$, then $$h((\alpha \vee \beta)) = T.$$

We third suppose that $h(\alpha) = F$ and $h(\beta) = F$. We need to prove that $$h((\alpha \vee \beta)) = F.$$
Hence $\langle\alpha,F\rangle$ and $\langle\beta,F\rangle$
both belong to $\mathscr{P}$, and therefore there exists a construction sequence $\langle\delta_1,\delta_2,...,\delta_m\rangle$ for $\alpha$ and its associated sequence of truth values $\langle\tau_1,\tau_2,...,\tau_m\rangle$ with the properties $\delta_m = \alpha$ and $\tau_m = F$. And there exists another construction sequence $\langle\varepsilon_1,\varepsilon_2, ..., \varepsilon_n\rangle$ for $\beta$ and its associated sequence of truth values $\langle\tau_1^{'},\tau_2^{'}, ..., \tau_n^{'}\rangle$ with the properties $\varepsilon_n = \beta$ and $\tau_n^{'} = F$. Hence $\langle\delta_1,\delta_2,...,\delta_m,\varepsilon_1,\varepsilon_2, ..., \varepsilon_n, (\alpha \vee \beta) \rangle$ is a construction sequence for $(\alpha \vee \beta)$. Because the condition 3. holds for the associated truth value sequences, the associated sequence of truth values for $\langle\delta_1,\delta_2,...,\delta_m,\varepsilon_1,\varepsilon_2, ..., \varepsilon_n, (\alpha \vee \beta) \rangle$ is $\langle\tau_1,\tau_2,...,F,\tau_1^{'},\tau_2^{'}, ..., F, F\rangle $, since $\tau_m = F$ and $\tau_n^{'} = F$. Hence, by the definition of $\mathscr{P}$, we obtain $\langle(\alpha \vee \beta), F\rangle \in \mathscr{P}$. Hence $h((\alpha \vee \beta)) = F$, as desired.

Therefore, Condition 3. holds for $h$. 

%leftarrow  leftrightarrow

\underline{Condition 4.} Suppose that $\alpha$ and $\beta$ both belong to $\overline{S}$. We need to show that 
$$h((\alpha \rightarrow \beta)) =
\setlength\arraycolsep{1pt}
\left\{
\begin{array}{rcrcrc@{\qquad}l}
F &\mbox{ if } h(\alpha) = T \mbox{ and } h(\beta) = F\\ 
T & \mbox{ otherwise }.\\ 
\end{array}
\right.
$$

We first suppose that $h(\alpha) = T$ and $h(\beta) = F$. We need to prove that $$h((\alpha \rightarrow \beta)) = F.$$
Hence $\langle\alpha,T\rangle$ and $\langle\beta,F\rangle$
both belong to $\mathscr{P}$, and therefore there exists a construction sequence $\langle\delta_1,\delta_2,...,\delta_m\rangle$ for $\alpha$ and its associated sequence of truth values $\langle\tau_1,\tau_2,...,\tau_m\rangle$ with the properties $\delta_m = \alpha$ and $\tau_m = T$. And there exists another construction sequence $\langle\varepsilon_1,\varepsilon_2, ..., \varepsilon_n\rangle$ for $\beta$ and its associated sequence of truth values $\langle\tau_1^{'},\tau_2^{'}, ..., \tau_n^{'}\rangle$ with the properties $\varepsilon_n = \beta$ and $\tau_n^{'} = F$. Hence $\langle\delta_1,\delta_2,...,\delta_m,\varepsilon_1,\varepsilon_2, ..., \varepsilon_n, (\alpha \rightarrow \beta) \rangle$ is a construction sequence for $(\alpha \rightarrow \beta)$. Because the condition 4. holds for the associated truth value sequences, the associated sequence of truth values for $\langle\delta_1,\delta_2,...,\delta_m,\varepsilon_1,\varepsilon_2, ..., \varepsilon_n, (\alpha \rightarrow \beta) \rangle$ is $\langle\tau_1,\tau_2,...,T,\tau_1^{'},\tau_2^{'}, ..., F, F\rangle $, since $\tau_m = T$ and $\tau_n^{'} = F$. Hence, by the definition of $\mathscr{P}$, we obtain $\langle(\alpha \rightarrow \beta), F\rangle \in \mathscr{P}$. Hence $h((\alpha \rightarrow \beta)) = F$, as desired.

We second suppose that $h(\alpha) = T$ and $h(\beta) = T$. We need to prove that $$h((\alpha \rightarrow \beta)) = T.$$
Hence $\langle\alpha,T\rangle$ and $\langle\beta,T\rangle$
both belong to $\mathscr{P}$, and therefore there exists a construction sequence $\langle\delta_1,\delta_2,...,\delta_m\rangle$ for $\alpha$ and its associated sequence of truth values $\langle\tau_1,\tau_2,...,\tau_m\rangle$ with the properties $\delta_m = \alpha$ and $\tau_m = T$. And there exists another construction sequence $\langle\varepsilon_1,\varepsilon_2, ..., \varepsilon_n\rangle$ for $\beta$ and its associated sequence of truth values $\langle\tau_1^{'},\tau_2^{'}, ..., \tau_n^{'}\rangle$ with the properties $\varepsilon_n = \beta$ and $\tau_n^{'} = T$. Hence $\langle\delta_1,\delta_2,...,\delta_m,\varepsilon_1,\varepsilon_2, ..., \varepsilon_n, (\alpha \rightarrow \beta) \rangle$ is a construction sequence for $(\alpha \rightarrow \beta)$. Because the condition 4. holds for the associated truth value sequences, the associated sequence of truth values for $\langle\delta_1,\delta_2,...,\delta_m,\varepsilon_1,\varepsilon_2, ..., \varepsilon_n, (\alpha \rightarrow \beta) \rangle$ is $\langle\tau_1,\tau_2,...,T,\tau_1^{'},\tau_2^{'}, ..., T, T\rangle $, since $\tau_m = T$ and $\tau_n^{'} = T$. Hence, by the definition of $\mathscr{P}$, we obtain $\langle(\alpha \rightarrow \beta), T\rangle \in \mathscr{P}$. Hence $h((\alpha \rightarrow \beta)) = T$, as desired.

We third suppose that $h(\alpha) = F$ and $h(\beta) = T$. We need to prove that $$h((\alpha \rightarrow \beta)) = T.$$
Hence $\langle\alpha,F\rangle$ and $\langle\beta,T\rangle$
both belong to $\mathscr{P}$, and therefore there exists a construction sequence $\langle\delta_1,\delta_2,...,\delta_m\rangle$ for $\alpha$ and its associated sequence of truth values $\langle\tau_1,\tau_2,...,\tau_m\rangle$ with the properties $\delta_m = \alpha$ and $\tau_m = F$. And there exists another construction sequence $\langle\varepsilon_1,\varepsilon_2, ..., \varepsilon_n\rangle$ for $\beta$ and its associated sequence of truth values $\langle\tau_1^{'},\tau_2^{'}, ..., \tau_n^{'}\rangle$ with the properties $\varepsilon_n = \beta$ and $\tau_n^{'} = T$. Hence $\langle\delta_1,\delta_2,...,\delta_m,\varepsilon_1,\varepsilon_2, ..., \varepsilon_n, (\alpha \rightarrow \beta) \rangle$ is a construction sequence for $(\alpha \rightarrow \beta)$. Because the condition 4. holds for the associated truth value sequences, the associated sequence of truth values for $\langle\delta_1,\delta_2,...,\delta_m,\varepsilon_1,\varepsilon_2, ..., \varepsilon_n, (\alpha \rightarrow \beta) \rangle$ is $\langle\tau_1,\tau_2,...,F,\tau_1^{'},\tau_2^{'}, ..., T, T\rangle $, since $\tau_m = F$ and $\tau_n^{'} = T$. Hence, by the definition of $\mathscr{P}$, we obtain $\langle(\alpha \rightarrow \beta), T\rangle \in \mathscr{P}$. Hence $h((\alpha \rightarrow \beta)) = T$, as desired.

We fourth suppose that $h(\alpha) = F$ and $h(\beta) = F$. We need to prove that $$h((\alpha \rightarrow \beta)) = T.$$
Hence $\langle\alpha,F\rangle$ and $\langle\beta,F\rangle$
both belong to $\mathscr{P}$, and therefore there exists a construction sequence $\langle\delta_1,\delta_2,...,\delta_m\rangle$ for $\alpha$ and its associated sequence of truth values $\langle\tau_1,\tau_2,...,\tau_m\rangle$ with the properties $\delta_m = \alpha$ and $\tau_m = F$. And there exists another construction sequence $\langle\varepsilon_1,\varepsilon_2, ..., \varepsilon_n\rangle$ for $\beta$ and its associated sequence of truth values $\langle\tau_1^{'},\tau_2^{'}, ..., \tau_n^{'}\rangle$ with the properties $\varepsilon_n = \beta$ and $\tau_n^{'} = F$. Hence $\langle\delta_1,\delta_2,...,\delta_m,\varepsilon_1,\varepsilon_2, ..., \varepsilon_n, (\alpha \rightarrow \beta) \rangle$ is a construction sequence for $(\alpha \rightarrow \beta)$. Because the condition 4. holds for the associated truth value sequences, the associated sequence of truth values for $\langle\delta_1,\delta_2,...,\delta_m,\varepsilon_1,\varepsilon_2, ..., \varepsilon_n, (\alpha \rightarrow \beta) \rangle$ is $\langle\tau_1,\tau_2,...,F,\tau_1^{'},\tau_2^{'}, ..., F, T\rangle $, since $\tau_m = F$ and $\tau_n^{'} = F$. Hence, by the definition of $\mathscr{P}$, we obtain $\langle(\alpha \rightarrow \beta), T\rangle \in \mathscr{P}$. Hence $h((\alpha \rightarrow \beta)) = T$, as desired.

Therefore, Condition 4. holds for $h$. 

\underline{Condition 5.} Suppose that $\alpha$ and $\beta$ both belong to $\overline{S}$. We need to show that 
$$h((\alpha \leftrightarrow \beta)) =
\setlength\arraycolsep{1pt}
\left\{
\begin{array}{rcrcrc@{\qquad}l}
T &\mbox{ if } h(\alpha) = h(\beta)\\ 
F & \mbox{ otherwise }.\\ 
\end{array}
\right.
$$

We first suppose that $h(\alpha) = h(\beta) $. We need to prove that $$h((\alpha \leftrightarrow \beta)) = T.$$ We consider two cases. 

Case1. Assume $h(\alpha) = h(\beta) = T$. 
Hence $\langle\alpha,T\rangle$ and $\langle\beta,T\rangle$
both belong to $\mathscr{P}$, and therefore there exists a construction sequence $\langle\delta_1,\delta_2,...,\delta_m\rangle$ for $\alpha$ and its associated sequence of truth values $\langle\tau_1,\tau_2,...,\tau_m\rangle$ with the properties $\delta_m = \alpha$ and $\tau_m = T$. And there exists another construction sequence $\langle\varepsilon_1,\varepsilon_2, ..., \varepsilon_n\rangle$ for $\beta$ and its associated sequence of truth values $\langle\tau_1^{'},\tau_2^{'}, ..., \tau_n^{'}\rangle$ with the properties $\varepsilon_n = \beta$ and $\tau_n^{'} = T$. Hence $\langle\delta_1,\delta_2,...,\delta_m,\varepsilon_1,\varepsilon_2, ..., \varepsilon_n, (\alpha \leftrightarrow \beta) \rangle$ is a construction sequence for $(\alpha \leftrightarrow \beta)$. Because the condition 5. holds for the associated truth value sequences, the associated sequence of truth values for $\langle\delta_1,\delta_2,...,\delta_m,\varepsilon_1,\varepsilon_2, ..., \varepsilon_n, (\alpha \leftrightarrow \beta) \rangle$
is $\langle\tau_1,\tau_2,...,T,\tau_1^{'},\tau_2^{'}, ..., T, T\rangle $, since $\tau_m = T$ and $\tau_n^{'} = T$. Hence, by the definition of $\mathscr{P}$, we obtain $\langle(\alpha \leftrightarrow \beta), T\rangle \in \mathscr{P}$. Hence $h((\alpha \leftrightarrow \beta)) = T$, as desired.

Case2. Assume $h(\alpha) = h(\beta) = F$. 
Hence $\langle\alpha,F\rangle$ and $\langle\beta,F\rangle$
both belong to $\mathscr{P}$, and therefore there exists a construction sequence $\langle\delta_1,\delta_2,...,\delta_m\rangle$ for $\alpha$ and its associated sequence of truth values $\langle\tau_1,\tau_2,...,\tau_m\rangle$ with the properties $\delta_m = \alpha$ and $\tau_m = F$. And there exists another construction sequence $\langle\varepsilon_1,\varepsilon_2, ..., \varepsilon_n\rangle$ for $\beta$ and its associated sequence of truth values $\langle\tau_1^{'},\tau_2^{'}, ..., \tau_n^{'}\rangle$ with the properties $\varepsilon_n = \beta$ and $\tau_n^{'} = F$. Hence $\langle\delta_1,\delta_2,...,\delta_m,\varepsilon_1,\varepsilon_2, ..., \varepsilon_n, (\alpha \leftrightarrow \beta) \rangle$ is a construction sequence for $(\alpha \leftrightarrow \beta)$. Because the condition 5. holds for the associated truth value sequences, the associated sequence of truth values for $\langle\delta_1,\delta_2,...,\delta_m,\varepsilon_1,\varepsilon_2, ..., \varepsilon_n, (\alpha \leftrightarrow \beta) \rangle$ is $\langle\tau_1,\tau_2,...,F,\tau_1^{'},\tau_2^{'}, ..., F, T\rangle $, since $\tau_m = F$ and $\tau_n^{'} = F$. Hence, by the definition of $\mathscr{P}$, we obtain $\langle(\alpha \leftrightarrow \beta), T\rangle \in \mathscr{P}$. Hence $h((\alpha \leftrightarrow \beta)) = T$, as desired.

We second suppose that $h(\alpha) \neq h(\beta) $. We need to prove that $$h((\alpha \leftrightarrow \beta)) = F.$$ We consider two cases. 

Case1. Assume $h(\alpha) = T$ and $h(\beta) = F$. 
Hence $\langle\alpha,T\rangle$ and $\langle\beta,F\rangle$
both belong to $\mathscr{P}$, and therefore there exists a construction sequence $\langle\delta_1,\delta_2,...,\delta_m\rangle$ for $\alpha$ and its associated sequence of truth values $\langle\tau_1,\tau_2,...,\tau_m\rangle$ with the properties $\delta_m = \alpha$ and $\tau_m = T$. And there exists another construction sequence $\langle\varepsilon_1,\varepsilon_2, ..., \varepsilon_n\rangle$ for $\beta$ and its associated sequence of truth values $\langle\tau_1^{'},\tau_2^{'}, ..., \tau_n^{'}\rangle$ with the properties $\varepsilon_n = \beta$ and $\tau_n^{'} = F$. Hence $\langle\delta_1,\delta_2,...,\delta_m,\varepsilon_1,\varepsilon_2, ..., \varepsilon_n, (\alpha \leftrightarrow \beta) \rangle$ is a construction sequence for $(\alpha \leftrightarrow \beta)$. Because the condition 5. holds for the associated truth value sequences, the associated sequence of truth values for $\langle\delta_1,\delta_2,...,\delta_m,\varepsilon_1,\varepsilon_2, ..., \varepsilon_n, (\alpha \leftrightarrow \beta) \rangle$ is $\langle\tau_1,\tau_2,...,T,\tau_1^{'},\tau_2^{'}, ..., F, F\rangle $, since $\tau_m = T$ and $\tau_n^{'} = F$. Hence, by the definition of $\mathscr{P}$, we obtain $\langle(\alpha \leftrightarrow \beta), F\rangle \in \mathscr{P}$. Hence $h((\alpha \leftrightarrow \beta)) = F$, as desired.

Case2. Assume $h(\alpha) = F$ and $h(\beta) = T$. We similarly prove $h((\alpha \leftrightarrow \beta)) = F$. 

Therefore, Condition 5. holds for $h$. 
\end{proof}

We have thus completed the proof that a function $\overline{v}$ satisfying the Conditions 0.-5., in the Theorem 12A, exists. It remains to prove the uniqueness part of the theorem. 

Let $S$ be any set of sentence symbols, and assume that $v : S \rightarrow \{F, T \}$ is a truth assignment. Show there is at most one extension $\overline{v}$ meeting Conditions 0.–5.. 

\begin{proof}
We prove the uniqueness part of the theorem by induction principle. Suppose $v$ is a truth assignment and $\overline{v_1},\overline{v_2}$ are two functions satisfying the conditions for $\overline{v}$. We prove by induction principle that $\overline{v_1}(\alpha) = \overline{v_2}(\alpha)$ for all wffs $\alpha$.

Base Case: If $\alpha$ is a sentence symbol belonging to $S$, then by Conditions 0., $\overline{v_1}(\alpha) = v(\alpha) = \overline{v_2}(\alpha)$. 

Induction Case: Suppose $\overline{v_1}(\alpha) = \overline{v_2}(\alpha)$. We need to show that $\overline{v_1}((\neg \alpha)) = \overline{v_2}((\neg \alpha))$. Assume $\overline{v_1}(\alpha) = T$. Then $\overline{v_1}((\neg \alpha)) = F$ by Condition 1.. Since $\overline{v_1}(\alpha) = \overline{v_2}(\alpha)$, then $\overline{v_2}(\alpha) = T$. Then $\overline{v_2}((\neg \alpha)) = F$, by Condition 1.. Thus, $\overline{v_1}((\neg \alpha)) = F = \overline{v_2}((\neg \alpha))$. The case of $\overline{v_1}(\alpha) = F$ is similar. Therefore, $\overline{v_1}((\neg \alpha)) = \overline{v_2}((\neg \alpha))$ has been proved.

Suppose $\overline{v_1}(\alpha) = \overline{v_2}(\alpha)$ and $\overline{v_1}(\beta) = \overline{v_2}(\beta)$. 

First, we need to show $\overline{v_1}((\alpha \wedge \beta)) = \overline{v_2}((\alpha \wedge \beta))$. Assume $\overline{v_1}(\alpha) = T$ and $\overline{v_1}(\beta) = T$. Then $\overline{v_1}((\alpha \wedge \beta)) = T$, by Condition 2.. Since $\overline{v_1}(\alpha) = \overline{v_2}(\alpha)$ and $\overline{v_1}(\beta) = \overline{v_2}(\beta)$, then $\overline{v_2}(\alpha) = T$ and $\overline{v_2}(\beta) = T$. Then $\overline{v_2}((\alpha \wedge \beta)) = T$, by Condition 2.. Thus, $\overline{v_1}((\alpha \wedge \beta)) = T = \overline{v_2}((\alpha \wedge \beta)) $. Therefore, $\overline{v_1}((\alpha \wedge \beta)) = \overline{v_2}((\alpha \wedge \beta))$. 

Next assume that not both $\overline{v_1}(\alpha) = T$ and $\overline{v_1}(\beta) = T$. Hence by Condition 2. $\overline{v_1}((\alpha \wedge \beta)) = F$. Also $\overline{v_2}(\alpha) = \overline{v_1}(\alpha)$ and
$\overline{v_2}(\beta) = \overline{v_1}(\beta)$, hence not both $\overline{v_2}(\alpha) = T$ and $\overline{v_2}(\beta) = T$. Hence, by Condition 2., $\overline{v_2}((\alpha \wedge \beta)) = F$. Thus, $\overline{v_1}((\alpha \wedge \beta)) = F = \overline{v_2}((\alpha \wedge \beta)) $. Therefore, $\overline{v_1}((\alpha \wedge \beta)) = \overline{v_2}((\alpha \wedge \beta))$. 

Second, we consider the case $\vee$.  
Assume $\overline{v_1}(\alpha) = T$ or $\overline{v_1}(\beta) = T$ (or both). Then $\overline{v_1}((\alpha \vee \beta)) = T$, by Condition 3.. Since $\overline{v_1}(\alpha) = \overline{v_2}(\alpha)$ and $\overline{v_1}(\beta) = \overline{v_2}(\beta)$, then $\overline{v_2}(\alpha) = T$ or $\overline{v_2}(\beta) = T$ (or both). Then $\overline{v_2}((\alpha \vee \beta)) = T$, by Condition 3.. Thus, $\overline{v_1}((\alpha \vee \beta)) = T = \overline{v_2}((\alpha \vee \beta)) $. Therefore, $\overline{v_1}((\alpha \vee \beta)) = \overline{v_2}((\alpha \vee \beta))$.

Next assume $\overline{v_1}(\alpha) = F$ and $\overline{v_1}(\beta) = F$. Then $\overline{v_1}((\alpha \vee \beta)) = F$, by Condition 3.. Since $\overline{v_1}(\alpha) = \overline{v_2}(\alpha)$ and $\overline{v_1}(\beta) = \overline{v_2}(\beta)$, then $\overline{v_2}(\alpha) = F$ and $\overline{v_2}(\beta) = F$. Then $\overline{v_2}((\alpha \vee \beta)) = F$, by Condition 3.. Thus, $\overline{v_1}((\alpha \vee \beta)) = F = \overline{v_2}((\alpha \vee \beta)) $. Therefore, $\overline{v_1}((\alpha \vee \beta)) = \overline{v_2}((\alpha \vee \beta))$.

Third, we consider the case $\rightarrow$. Assume $\overline{v_1}(\alpha) = T$ and $\overline{v_1}(\beta) = F$. Then $\overline{v_1}((\alpha \rightarrow \beta)) = F$, by Condition 4.. Since $\overline{v_1}(\alpha) = \overline{v_2}(\alpha)$ and $\overline{v_1}(\beta) = \overline{v_2}(\beta)$, then $\overline{v_2}(\alpha) = T$ and $\overline{v_2}(\beta) = F$. Then $\overline{v_2}((\alpha \rightarrow \beta)) = F$, by Condition 4.. Thus, $\overline{v_1}((\alpha \rightarrow \beta)) = F = \overline{v_2}((\alpha \rightarrow \beta)) $. Therefore, $\overline{v_1}((\alpha \rightarrow \beta)) = \overline{v_2}((\alpha \rightarrow \beta))$.

Next we consider the three remaining assumptions of $\rightarrow$. Those are $\overline{v_1}(\alpha) = F$ and $\overline{v_1}(\beta) = T$, $\overline{v_1}(\alpha) = \overline{v_1}(\beta) = T$, and $\overline{v_1}(\alpha) = \overline{v_1}(\beta) = F$. By Condition 4., $\overline{v_1}((\alpha \rightarrow \beta)) = T$ satisfies for all the three assumptions. Since $\overline{v_1}(\alpha) = \overline{v_2}(\alpha)$ and $\overline{v_1}(\beta) = \overline{v_2}(\beta)$, then $\overline{v_2}(\alpha)$ and $\overline{v_2}(\beta)$ satisfies one the three assumptions, those are $\overline{v_2}(\alpha) = F$ and $\overline{v_2}(\beta) = T$, $\overline{v_2}(\alpha) = \overline{v_2}(\beta) = T$, or $\overline{v_2}(\alpha) = \overline{v_2}(\beta) = F$. Then $\overline{v_2}((\alpha \rightarrow \beta)) = T$, by Condition 4.. Thus, $\overline{v_1}((\alpha \rightarrow \beta)) = T = \overline{v_2}((\alpha \rightarrow \beta)) $. Therefore, $\overline{v_1}((\alpha \rightarrow \beta)) = \overline{v_2}((\alpha \rightarrow \beta))$.

Fourth, we consider the case $\leftrightarrow$. Assume $\overline{v_1}(\alpha) = \overline{v_1}(\beta)$. Then $\overline{v_1}((\alpha \leftrightarrow \beta)) = T$, by Condition 5.. Since $\overline{v_1}(\alpha) = \overline{v_2}(\alpha)$ and $\overline{v_1}(\beta) = \overline{v_2}(\beta)$, then $\overline{v_2}(\alpha) = \overline{v_2}(\beta)$. Then $\overline{v_2}((\alpha \leftrightarrow \beta)) = T$, by Condition 5.. Thus, $\overline{v_1}((\alpha \leftrightarrow \beta)) = T = \overline{v_2}((\alpha \leftrightarrow \beta)) $. Therefore, $\overline{v_1}((\alpha \leftrightarrow \beta)) = \overline{v_2}((\alpha \leftrightarrow \beta))$. 

Next assume $\overline{v_1}(\alpha) \neq \overline{v_1}(\beta)$. Then $\overline{v_1}((\alpha \leftrightarrow \beta)) = F$, by Condition 5.. Since $\overline{v_1}(\alpha) = \overline{v_2}(\alpha)$ and $\overline{v_1}(\beta) = \overline{v_2}(\beta)$, then $\overline{v_2}(\alpha) \neq \overline{v_2}(\beta)$. Then $\overline{v_2}((\alpha \leftrightarrow \beta)) = F$, by Condition 5.. Thus, $\overline{v_1}((\alpha \leftrightarrow \beta)) = F = \overline{v_2}((\alpha \leftrightarrow \beta)) $. Therefore, $\overline{v_1}((\alpha \leftrightarrow \beta)) = \overline{v_2}((\alpha \leftrightarrow \beta))$. 

Thus, we proved that there is at most one extension $v$ meeting conditions 0.–5.. Therefore, we proved the uniqueness part of the theorem.

\end{proof}

\section{Recursion}
\label{Mar. 13th}
%We return now to the more abstract case. There is a set $U$ (such as the set of all expressions), a subset $B$ of $U$ (such as the set of sentence symbols), and two functions $f$ and $g$, where $$f:U \times U \rightarrow U \mbox{ and } g: U \rightarrow U.$$

%{\bf $C$ is the set generated from $B$ by $f$ and $g$.} 

We want to consider a general situation where we are given a set $U$, and a family $\mathscr{F}$ of functions such that for every $f\in \mathscr{F}$, the domain of $f$ is $U^n$ for some $n \succeq 0$, the co-domain of $f$ is $U$, i.e. $f: U^n \rightarrow U$. Such a function $f$ is called n-ary. For the sake of clarity, when $n=0$, then we consider $f$ to be an element of $U$, i.e. $f \in U$ is a constant. 

We are furthermore given a subset $B$ of $U$. We now make the following Definition. 

{\bf Definition 1.} A subset $C \subseteq U$ is generated from $B$ by $\mathscr{F}$ if the following three conditions hold: 

(a) $B \subseteq C$,

(b) For every n-ary function $f\in \mathscr{F}$, $f(C^n) \subseteq C$, and 

(c) If $C'$ is a proper subset of $C$ such that $B \subseteq C'$, then for some natural number $n$, and an n-ary $f \in \mathscr{F}$, $f((C')^n) \nsubseteq C'$. 

{\bf Comment.} Any subset $C$ of $U$, for which the condition (b) holds, is called $\mathscr{F}$-closed. Thus the conditions (b) and (c) together say that C is $\mathscr{F}$-closed, but no proper subset of $C$, which contains B, is $\mathscr{F}$-closed. Furthermore, $C$ is generated from $B$ by $\mathscr{F}$ if $B \subseteq C$, $C$ is $\mathscr{F}$-closed, and no proper subset of $C$, containing $B$, is $\mathscr{F}$-closed. 

{\bf Theorem.} With the above notation and Definitions, there exists a unique subset $C$, of $U$, generated from $B$ by $\mathscr{F}$.

{\bf Definition 2.} With the notation introduced above, we say that a subset $C$, of $U$, is {\it freely generated} from $B$ by $\mathscr{F}$ if $C$ is generated from $B$ by $\mathscr{F}$ (see Definition 1.), and the following three conditions holds:

\begin{itemize}
\item \texttt{F1.} Every $f$ in $\mathscr{F}$ is one-to-one on $C$;
\item \texttt{F2.} For every $f$ in $\mathscr{F}$, the range of $f$ on $C$ is disjoint from $B$, i.e. $f(C)\cap B = \emptyset$; and 
\item \texttt{F3.} For every $f_1$ and $f_2$ in $\mathscr{F}$, if $f_1 \neq f_2$, then the ranges of $f_1$ and $f_2$ are disjoint on $C$, i.e. $f_1(C)\cap f_2(C) = \emptyset$. 
\end{itemize}

We frequently encounter the following situation. We are given a function $h: B \rightarrow V$, where $V$ is some set of ``values'', and we want to extend, $h$ to a function $\overline{h}: C \rightarrow V$, assigning ``values'' to all elements of $C$, according to some ``rules'' associated with the functions in $\mathscr{F}$: Such an extension $\overline{h}$ exists when $C$ is {\it freely generated} from $B$ by $\mathscr{F}$. The following Recursion Theorem contains a rigorous formulation of the preceding informal comments. 

In addition to the function $h: B \rightarrow V$, we are also given, for each $f \in \mathscr{F}$, a corresponding ``rule'' $\tilde{f}$, which has the property that, if $f$ is n-ary, then $\tilde{f}: V^n \rightarrow V$. 

{\bf Recursion Theorem.} Suppose that $C$ is {\it freely generated} from $B$ by $\mathscr{F}$, $h: B \rightarrow V$, and each $f$ in $\mathscr{F}$ has an associated $\tilde{f}$ as described above. 

Then there is a unique function $\overline{h} : C \rightarrow V$ satisfying the following conditions: 

(i) For $x$ in $B$, $\overline{h}(x) = h(x)$; 

(ii) For every $f \in \mathscr{F}$, if $f$ is n-ary, then for every $x_1, x_2, ..., x_n \in C$, $\overline{h}(f(x_1, x_2, ..., x_n)) = \tilde{f}(\overline{h}(x_1), \overline{h}(x_2), ..., \overline{h}(x_n))$. 

We will now show how the Truth Assignment Theorem (Theorem 12 A, p.23 in \cite{enderton2001mathematical}) can be obtained as a consequence of the Recursion Theorem. We begin by specifying the parameters for the desired application of the Recursion Theorem. 

We will set $U$ to be the set of all LSL-expressions, and $\mathscr{F} = \{\varepsilon_\neg,  \varepsilon_\wedge,\varepsilon_\vee, \varepsilon_\rightarrow, \varepsilon_\leftrightarrow \}$ be the set of the formula-building operations, as explained at the beginning of the paper, on p.2. Thus $\varepsilon_\neg: U\rightarrow U$ is unary, and the remaining four $\varepsilon_\wedge, \varepsilon_\vee, \varepsilon_\rightarrow, \varepsilon_\leftrightarrow$ mapping $U^2 \rightarrow U$, are binary. 

We also set $B=S$, where $S$ is an arbitrary set of sentence symbols. 

Furthermore we set $V = \{T,F\}$ to be the set of truth values and we define $\tilde{\varepsilon_\neg}: V\rightarrow V $, and $ \tilde{\varepsilon_\wedge}, \tilde{\varepsilon_\vee}, \tilde{\varepsilon_\rightarrow}, \tilde{\varepsilon_\leftrightarrow} : V^2\rightarrow V$ in accord with the conditions 1.-5. on pages 20, 21 in \cite{enderton2001mathematical}, i.e. $\tilde{\varepsilon_\neg}(T) = F, \tilde{\varepsilon_\neg} (F)= T; \tilde{\varepsilon_\wedge}(T,T) = T, \tilde{\varepsilon_\wedge}(T,F) = \tilde{\varepsilon_\wedge}(F,T) = \tilde{\varepsilon_\wedge}(F,F) = F; \tilde{\varepsilon_\vee}(T,T) = \tilde{\varepsilon_\vee}(T,F) = \tilde{\varepsilon_\vee}(F,T) = T, \tilde{\varepsilon_\vee}(F,F) = F; \tilde{\varepsilon_\rightarrow}(T,T) = \tilde{\varepsilon_\rightarrow}(F,T) = \tilde{\varepsilon_\rightarrow}(F,F) = T, \tilde{\varepsilon_\rightarrow}(T,F) = F; \tilde{\varepsilon_\leftrightarrow}(T,T) = \tilde{\varepsilon_\leftrightarrow}(F,F) = T, \tilde{\varepsilon_\leftrightarrow}(T,F) = \tilde{\varepsilon_\leftrightarrow}(F,T) = F.$

We now observe that the set $\overline{S}$, which is the set of all wffs $\alpha$, such that every sentence symbol occurring in $\alpha$ belongs to $S$, is a freely generated subset of the set $U$ of all LSL expressions, with the set $B$ of the generating elements being $S$, and the set of generating functions $\mathscr{F}$ being $\{ \varepsilon_\neg, \varepsilon_\wedge, \varepsilon_\vee, \varepsilon_\rightarrow, \varepsilon_\leftrightarrow \}$. This fact is simply the Unique Readability Theorem for the wffs of LSL. In addition to the class notes, this Theorem is stated on page 40 in \cite{enderton2001mathematical}, and we also rely on this theorem in the proof of the Main Lemma earlier in the paper. Thus to obtain the Truth Assignment Theorem as an application of the Recursion Theorem, we will set $C= \overline{S}$. We finally set $h: B \rightarrow V$ to be an arbitrary truth assignment $v:S \rightarrow \{T, F\}$. It is then clear from the preceding discussion that all assumptions of the Recursion Theorem hold, and thus there exists a unique function $\overline{v} = \overline{h}: \overline{S} \rightarrow \{T, F\}$ satisfying the conditions (i) and (ii) of the Recursion Theorem. 

It thus remains to verify that the condition function $\overline{v} : \overline{S} \rightarrow \{T,F\}$ we obtained. Similarly as in our first prove of the Theorem 12A, it remains to verify that the conditions 0.-5. hold for the function $\overline{v}:\overline{S}\rightarrow\{T,F\}$ we just obtained. 

We now proceed to do so. 

\underline{Condition 0.} Suppose that $\alpha$ is a sentence symbol belonging to $S$, i.e. $\alpha = A_n \in S$. But in our application of the Recursion Theorem, we have set $B = S$, and thus by the condition (i) in the statement of the Recursion Theorem, we obtain $\overline{h}(\alpha) = h(\alpha)$, i.e. $\overline{v}(\alpha) = v(\alpha)$, having set $h = v$. Hence the Condition 0. holds. 

\underline{Condition 1.} Let $\alpha \in \overline{S}$. We need to show that 

$$\overline{v}(( \neg \alpha)) =
\setlength\arraycolsep{1pt}
\left\{
\begin{array}{rcrcrc@{\qquad}l}
T &\mbox{ if } \overline{v}(\alpha) = F\\ 
F &\mbox{ if } \overline{v}(\alpha) = T.\\ 
\end{array}
\right.
$$

In our application of the Recursion Theorem we have set $C= \overline{S}$, and we have $\varepsilon_\neg \in \mathscr{F}$, and the condition (b) of Definition 1 holds. Hence, $\epsilon_\neg (\alpha) = (\neg \alpha) \in C = \overline{S}$. Thus by the condition (ii) of the Recursion Theorem, we obtain $\overline{h}(\varepsilon_\neg (\alpha)) = \tilde{\varepsilon_\neg}(\overline{h}(\alpha))$, i.e. 

$$\overline{v}(( \neg \alpha)) =\tilde{\varepsilon_\neg}(\overline{v}(\alpha))=
\setlength\arraycolsep{1pt}
\left\{
\begin{array}{rcrcrc@{\qquad}l}
T &\mbox{ if } \overline{v}(\alpha) = F\\ 
F &\mbox{ if } \overline{v}(\alpha) = T\\ 
\end{array}
\right.
$$ in accord with the specification of $\tilde{\varepsilon_\neg}$. Thus the Condition 1. holds.

\underline{Condition 2.} Let $\alpha$ and $\beta$ both belong to $\overline{S}$. We need to show that 
$$\overline{v}((\alpha \wedge \beta)) =
\setlength\arraycolsep{1pt}
\left\{
\begin{array}{rcrcrc@{\qquad}l}
T &\mbox{ if } \overline{v}(\alpha) = \overline{v}(\beta) = T\\ 
F & \mbox{ otherwise. }\\ 
\end{array}
\right.
$$

In our application of the Recursion Theorem we have set $C= \overline{S}$, and we have $\varepsilon_\wedge \in \mathscr{F}$, and the condition (b) of Definition 1 holds. Hence, $\varepsilon_\wedge (\alpha, \beta) = (\alpha \wedge \beta) \in C = \overline{S}$. Thus by the condition (ii) of the Recursion Theorem, we obtain $\overline{h}(\varepsilon_\wedge (\alpha, \beta)) = \tilde{\varepsilon_\wedge}(\overline{h}(\alpha), \overline{h}(\beta))$, i.e. $$\overline{v}((\alpha \wedge \beta)) = \tilde{\varepsilon_\wedge}(\overline{v}(\alpha), \overline{v}(\beta))=
\setlength\arraycolsep{1pt}
\left\{
\begin{array}{rcrcrc@{\qquad}l}
T &\mbox{ if } \overline{v}(\alpha) = \overline{v}(\beta) = T\\ 
F & \mbox{ otherwise }\\ 
\end{array}
\right.
$$ in accord with the specification of $\tilde{\varepsilon_\wedge}$. Thus the Condition 2. holds.

\underline{Condition 3.} Let $\alpha$ and $\beta$ both belong to $\overline{S}$. We need to show that 
$$\overline{v}((\alpha \vee \beta)) =
\setlength\arraycolsep{1pt}
\left\{
\begin{array}{rcrcrc@{\qquad}l}
T &\mbox{ if } \overline{v}(\alpha) = T \mbox{ or } \overline{v}(\beta) = T \mbox{ (or both) }\\ 
F & \mbox{ otherwise. }\\ 
\end{array}
\right.
$$

In our application of the Recursion Theorem we have set $C= \overline{S}$, and we have $\varepsilon_\vee \in \mathscr{F}$, and the condition (b) of Definition 1 holds. Hence, $\varepsilon_\vee (\alpha, \beta) = (\alpha \vee \beta) \in C = \overline{S}$. Thus by the condition (ii) of the Recursion Theorem, we obtain $\overline{h}(\varepsilon_\vee (\alpha, \beta)) = \tilde{\varepsilon_\vee}(\overline{h}(\alpha), \overline{h}(\beta))$, i.e. $$\overline{v}((\alpha \vee \beta)) = \tilde{\varepsilon_\vee}(\overline{v}(\alpha), \overline{v}(\beta))=
\setlength\arraycolsep{1pt}
\left\{
\begin{array}{rcrcrc@{\qquad}l}
T &\mbox{ if } \overline{v}(\alpha) = T \mbox{ or } \overline{v}(\beta) = T \mbox{ (or both) }\\ 
F & \mbox{ otherwise }\\ 
\end{array}
\right.
$$ in accord with the specification of $\tilde{\varepsilon_\vee}$. Thus the Condition 3. holds.

\underline{Condition 4.} Let $\alpha$ and $\beta$ both belong to $\overline{S}$. We need to show that 
$$\overline{v}((\alpha \rightarrow \beta)) =
\setlength\arraycolsep{1pt}
\left\{
\begin{array}{rcrcrc@{\qquad}l}
F &\mbox{ if } \overline{v}(\alpha) = T \mbox{ and } \overline{v}(\beta) = F\\ 
T & \mbox{ otherwise. }\\ 
\end{array}
\right.
$$

In our application of the Recursion Theorem we have set $C= \overline{S}$, and we have $\varepsilon_\rightarrow \in \mathscr{F}$, and the condition (b) of Definition 1 holds. Hence, $\varepsilon_\rightarrow (\alpha, \beta) = (\alpha \rightarrow \beta) \in C = \overline{S}$. Thus by the condition (ii) of the Recursion Theorem, we obtain $\overline{h}(\varepsilon_\rightarrow (\alpha, \beta)) = \tilde{\varepsilon_\rightarrow}(\overline{h}(\alpha), \overline{h}(\beta))$, i.e. $$\overline{v}((\alpha \rightarrow \beta)) = \tilde{\varepsilon_\rightarrow}(\overline{v}(\alpha), \overline{v}(\beta))=
\setlength\arraycolsep{1pt}
\left\{
\begin{array}{rcrcrc@{\qquad}l}
F &\mbox{ if } \overline{v}(\alpha) = T \mbox{ and } \overline{v}(\beta) = F\\ 
T & \mbox{ otherwise }\\ 
\end{array}
\right.
$$ in accord with the specification of $\tilde{\varepsilon_\rightarrow}$. Thus the Condition 4. holds.

\underline{Condition 5.} Let $\alpha$ and $\beta$ both belong to $\overline{S}$. We need to show that 
$$\overline{v}((\alpha \leftrightarrow \beta)) =
\setlength\arraycolsep{1pt}
\left\{
\begin{array}{rcrcrc@{\qquad}l}
T &\mbox{ if } \overline{v}(\alpha) = \overline{v}(\beta)\\ 
F & \mbox{ otherwise. }\\ 
\end{array}
\right.
$$

In our application of the Recursion Theorem we have set $C= \overline{S}$, and we have $\varepsilon_\leftrightarrow \in \mathscr{F}$, and the condition (b) of Definition 1 holds. Hence, $\varepsilon_\leftrightarrow (\alpha, \beta) = (\alpha \leftrightarrow \beta) \in C = \overline{S}$. Thus by the condition (ii) of the Recursion Theorem, we obtain $\overline{h}(\varepsilon_\leftrightarrow (\alpha, \beta)) = \tilde{\varepsilon_\leftrightarrow}(\overline{h}(\alpha), \overline{h}(\beta))$, i.e. $$\overline{v}((\alpha \leftrightarrow \beta)) = \tilde{\varepsilon_\leftrightarrow}(\overline{v}(\alpha), \overline{v}(\beta))=
\setlength\arraycolsep{1pt}
\left\{
\begin{array}{rcrcrc@{\qquad}l}
T &\mbox{ if } \overline{v}(\alpha) = \overline{v}(\beta)\\ 
F & \mbox{ otherwise }\\ 
\end{array}
\right.
$$ in accord with the specification of $\tilde{\varepsilon_\leftrightarrow}$. Thus the Condition 5. holds.

\section{Conclusion}
\label{conclusion_chapter}
In this paper, we have proved the Truth Assignment Theorem for the language of sentential logic. We have done this in two different ways. In the first approach, we have proved the theorem ``directly'', i.e. we have provided a complete proof with the exception that we are using the Unique Readability Theorem for the language of sentential logic without proof. In the second approach, we derive the Truth Assignment Theorem as a consequence of a general formulation of the Recursion Theorem. We again use the Unique Readability Theorem without proof. 

\clearpage

\addcontentsline{toc}{section}{References}

\printbibliography

\end{document}